\documentclass[12pt]{article}
\usepackage{amsmath, latexsym, amssymb}

\newcommand{\R}{\mathbb R}
\newcommand{\C}{\mathbb C}

\newcommand{\QED}{\hspace{.2in}\square\newline}
\newcommand{\qED}{\hspace{.2in}\boxminus\newline}

\newtheorem{theorem}{Theorem}[section]
\newtheorem{corollary}{Corollary}[section]

\newtheorem{claim}{Claim}[section]

\newtheorem{lemma}{Lemma}[section]

\begin{document}

\begin{center}
{\Large \textbf{Functional Mellin and Zeta Zeros}} \vskip 4em
{\large J. LaChapelle}\\
\vskip 2em
\end{center}

\begin{abstract}
A key theorem formulated in the context of functional Mellin transforms generalizes the important relationship $\exp\mathrm{tr} M=\det\exp M$. Along with the involution symmetry of the zeta function, the theorem suggests a strategy for tackling the Riemann hypothesis.
\end{abstract}

\section{Introduction}
 A functional integral scheme based on topological groups and Banach algebras was proposed in \cite{LA2} and subsequently used to define functional Mellin transforms \cite{LA1}. Such transforms are useful tools for probing Banach $\ast$-algebras. Specifically, they can be used to represent \emph{complex} powers, traces, and determinants of linear, bounded operators on some Hilbert space. They play a central role in a key theorem that extends the well-known relation $\exp\mathrm{tr} M=\det\exp M$ to general Banach $\ast$-algebras: Importantly, the extension holds for complex powers as well.

In particular, the  theorem can be applied to an operator on some Hilbert space whose spectrum is the natural numbers. Meanwhile, the trace of an appropriate functional Mellin transform associated with such an operator yields the Riemann zeta function, and the key theorem along with the symmetry of the zeta zeros under involution \emph{formally} implies the zeros must have real part equal to $1/2$. This of course is the Riemann hypothesis. Since the road to its proof is infamously treacherous, the proposed strategy employing the key theorem is submitted with due caution.\footnote{Admittedly, the strategy may be too simple to capture the full essence of the problem.} Even if the strategy is viable, we only claim consistency with the hypothesis --- not proof: After all, our motivation and presentation are at a formal level, say of mathematical physics, and not rigorous functional analysis.

\section{$\Re(\alpha_0)=1/2$}
As a brief preview; the idea is to construct an operator $A^{-\alpha}$ (with $\alpha\in\C$) acting on a suitable Hilbert space whose trace coincides with the Riemann zeta function and then to make use of the key theorem, which in this case essentially boils down to $e^{-\mathrm{tr}\, A^{-\alpha}}=\det e^{-A^{-\alpha}}$. Together with the involution symmetry $e^{-\mathrm{tr}\, A^{-\alpha_0}}=e^{-\mathrm{tr}\, A^{1-\alpha^\ast_0}}$ (and suitable conditions), we argue this implies $\Re(\alpha_0)=1/2$ where $\alpha_0$ is a zeta zero.

Most of the background and notation for the relevant aspects of functional Mellin will be assumed but can be found in \cite{LA1}.

\subsection{Eta and zeta as functional Mellin}
Let us quickly review the well-known Mellin representations of Dirichlet eta and Riemann zeta.

Consider a self-adjoint $A\in L_B(\mathcal{H})$  with $\sigma(A)=\mathbb{N}_+$. Let
$\{|j\rangle,\varepsilon_j\}$ with $j\in\{1,\ldots,\infty\}$ denote
the set of orthonormal eigenvectors and eigenvalues of
$A$. The Dirichlet eta
function associated with $A$ is given by
\begin{eqnarray}
\eta_{\mathrm{A}_{\Gamma}}(\alpha)
&=&\mathrm{tr}\int_{\R_+}-e^{-A(x+i\pi)}x^\alpha
\;d\nu(x_{\Gamma})\notag\\
&=&-\sum_{j=1}^\infty\int_{\R_+}\,e^{-\varepsilon_j (x+i\pi)+\alpha(\log x)}\langle j|j\rangle
\;d\nu(x_{\Gamma})\notag\\
&=&\int_0^\infty\frac{1}{e^{\,x}+1}\,x^{\alpha} \;d\nu(x_{\Gamma})
 \;,\;\;\;\;\;\alpha\in\langle0,\infty\rangle_{\Gamma}
\end{eqnarray}
where $\nu(x_{\Gamma}):=\nu(x)/\Gamma(\alpha)$ with $\nu(x)$ the normalized Haar measure on $\R_+$.

Likewise, for the same operator $A\in L_B(\mathcal{H})$ with $A=A^\ast$ and
$\sigma(A)=\mathbb{N}_+$, the Riemann zeta function associated with
$A$ can be defined by
\begin{eqnarray}
\zeta_{\mathrm{A}_\Gamma}(\alpha)
&=&\mathrm{tr}\int_{\R_+}e^{-Ax}x^{\alpha}
\;d\nu(x_{{\Gamma}})\notag\\
&=&\sum_{j=1}^\infty\int_{\R_+}\,e^{-\varepsilon_j x+\alpha(\log x)}\langle j|j\rangle
\;d\nu(x_{{\Gamma}})\notag\\
&=&\int_0^\infty\frac{1}{e^{\,x}-1}\,x^{\alpha} \;d\nu(x_{\Gamma})
\;,\;\;\;\;\;\alpha\in\langle1,\infty\rangle_{\Gamma}
\end{eqnarray}
where again $\nu(x_{\Gamma}):=\nu(x)/\Gamma(\alpha)$.

To make use of the theorem in the next section, we need a proper
 functional Mellin representation of Riemann zeta. Functional Mellin requires the data $\{G^\C,\mathfrak{C}^\ast,G^\C_\Lambda\}$. Let $G^\C$ be the space of continuous \emph{pointed} homomorphisms $\Omega \C^\times:=\mathrm{Hom}_C((\mathrm{S}^1,s_0),(\C^\times,z_0))$. When equipped with the usual topology, this is a topological abelian group with the uniform convergence topology on compact sets.  For the  $C^\ast$-algebra $\mathfrak{C}^\ast$, we take $\mathfrak{Z}$ which denotes the abelian group of integers $\mathbb{Z}$ over the complex numbers equipped with the usual ring structure of the integers. In this application, a convenient choice for $G^\C_\Lambda$ is the single locally compact topological group  $G^\C_\lambda\cong\mathcal{C}$ where $\mathcal{C}\cong S^1$ is a smooth \emph{clockwise} contour in $\C^\times$. Next, we need a representation $\rho:\mathcal{C}\rightarrow\mathfrak{Z}$. Since $g$ is a loop labelled by its degree, the appropriate representation to choose is $\rho(g_\lambda)=z Id$ where $z\in\C^\times$ and $Id$ is the identity on algebra $\mathfrak{Z}$.

Let $\mathrm{A}:\Omega \C^\times\rightarrow \mathfrak{Z}$ be the functional $g\mapsto \mathrm{A}(g):=\mathrm{deg}(g)\,g$ with $g\in\Omega \C^\times$. Restricting $\mathrm{A}$ to $G^\C_\lambda$ induces a function $A:\mathcal{C}\rightarrow \mathfrak{Z}$ given by $g_\lambda\mapsto A\rho(g_\lambda)={deg}(g_\lambda)\rho(g_\lambda)$. Now we need a representation $\pi:\mathfrak{Z}\rightarrow L_B(\mathcal{Z})$ where $\mathcal{Z}$ denotes the module $\mathbb{Z}$ over the real numbers. Let $|n\rangle$ be an orthonormal basis in $\mathcal{Z}$ labeling homotopy equivalence classes of pointed loops in $\C^\times$. Then $\pi'(\mathrm{A}(g))|n\rangle=\pi'(A\rho(g_\lambda))|n\rangle=nz|n\rangle$ with $n\in\mathbb{N}$. For Riemann zeta, we need to calculate the functional trace of $\mathrm{A}$ which means the zero spectrum of $\mathrm{A}(g)$ must be excluded.\cite{LA1} Thus we have $\langle n'|\pi(e^{-\mathrm{A}(\rho(g_\lambda))})|n\rangle=\langle n'|e^{-z\pi'( A)}|n\rangle=e^{-nz}\delta_{n' n}$ now with $n\in\mathbb{N}_+$.

 In a diagonal basis then (and branch cut for $\log$ along the positive reals),
\begin{eqnarray}
\mathcal{M}_{\mathcal{C}}\left[\mathrm{Tr}\,\mathrm{E}^{-\mathrm{A}};\alpha\right]
&=&\int_{\Omega \C^\times}\mathrm{tr}\,(e^{-\mathrm{A}(g)}\rho(g^\alpha))\;\mathcal{D}_{\mathcal{C}} g\notag\\
&\stackrel{\lambda}{\rightarrow}&\int_{\mathcal{C}}\mathrm{tr}\,(e^{-\mathrm{A}(g_\mathcal{C})}
\rho(g_\mathcal{C}^\alpha))\;d\nu(g_\mathcal{C})\notag\\
&=&\int_{\mathcal{C}}\sum_{n=1}^\infty(e^{-nz})\,z^{\alpha}\;d\nu(z_\mathcal{C})\notag\\
&=&\int_{\R_+}\frac{1}{e^z-1}\,r^{\alpha}\;d\nu(r_{\mathcal{C}}),
\;\;\;\;\;\alpha\in\langle0,\infty\rangle\backslash\{1\}\notag\\
&=&\zeta_{\mathrm{A}_\mathcal{C}}(\alpha)
\end{eqnarray}
where $d\nu(r_{\mathcal{C}})
:=\frac{\pi\csc(\pi\alpha)}{2\pi\imath\,\Gamma(\alpha)}\frac{dr}{r}$ and
$\mathcal{C}$ starts at $+\infty$ just below the real axis, passes
around the origin clockwise, and then continues back to
$+\infty$ above the real axis.\footnote{This assumes $\arg z=\pm\pi$ below(respectively above) the real axis. Of course
$\pi\csc(\pi\alpha)/\Gamma(\alpha)$ can be analytically continued to
the left-half plane thus obtaining a meromorphic representation of
Riemann zeta.}

\subsection{Key theorem}
\begin{theorem}\label{key theorem}
Suppose $\mathrm{A}\in \mathrm{Mor}_C(G^\C,\mathfrak{C}^\ast)$  by $g\mapsto A (g)$ is trace class  such that $0\notin\sigma(A)$, and let $\mathrm{E}^{-\mathrm{A}}\in \mathrm{Mor}_C(G^\C,\mathfrak{C}^\ast)$  by $g\mapsto e^{-A}(g)$ with $e^{-A}\in \mathfrak{C}^\ast$. Assume that $\mathrm{Tr}\, \mathrm{E}^{-\mathrm{A}}$ and $\mathrm{E}^{-\mathrm{Tr}\,\mathrm{A}}$ are Mellin integrable for a common domain
$\alpha\in\mathbb{S}_\lambda$. Then
\begin{equation}
e^{-\mathrm{tr}\,A_\lambda^{-\alpha}} =(\mathrm{Det}\,
\mathrm{E}^{-\mathrm{A}})_\lambda^{-\alpha}\;,\;\;\;\;\;\alpha\in\mathbb{S}_\lambda\;.
\end{equation}
\end{theorem}
Roughly stated, this theorem expresses the conditions under which the functional Mellin transform and exponential map commute. Depending on what is $A$, there may or may not exist a choice of $\lambda$ and an associated fundamental region $\mathbb{S}_\lambda$ where the equality holds.
\vskip 1em

\emph{Proof sketch}\footnote{The full proof including supporting propositions can be found in \cite{LA1}.}: First, recall that an
immediate consequence of the definitions and the relationship
between $(\mathrm{exp}\,\mathrm{tr})$ and
$(\det\,\mathrm{exp})$ is
\begin{eqnarray}
\mathcal{M}_{\lambda}
\left[\mathrm{Det}\,\mathrm{E}^{-\mathrm{A}};\alpha\right]
&=&\int_{G^\C}\det\left(e^{-A(g)
}g^\alpha\right)\,\mathcal{D}_\lambda g\notag\\
&=&\int_{G^\C} e^{-\mathrm{tr}\,A(g)}{\det} g^\alpha\,\mathcal{D}_\lambda g\notag\\
&=&\mathcal{M}_{\lambda}
\left[\mathrm{E}^{-\mathrm{Tr}\,\mathrm{A}};\alpha\right]
\end{eqnarray}
where the second line follows as soon as $A(g)$ is
trace class.

\begin{lemma}\label{key lemma}
Let $\mathrm{Tr}\, \mathrm{F}\in \mathrm{Mor}_C(G^\C,\C)$. Suppose the
functional Mellin transforms of $\,\mathrm{Tr}\, \mathrm{F}$ and
$\mathrm{E}^{-\mathrm{Tr}\, \mathrm{F}}$ exist for common $\alpha\in\mathbb{S}_\lambda$ for a given $\lambda$. Then
\begin{equation}
\mathcal{M}_{\lambda}\left[\mathrm{E}^{-\mathrm{Tr}\,
\mathrm{F}};\alpha\right] =e^{-\mathcal{M}_{\lambda}\left[\, \mathrm{Tr}\,
\mathrm{F};\alpha\right]}\,,\;\;\;\;\;\alpha\in\mathbb{S}_\lambda\;.
\end{equation}
\end{lemma}
\emph{proof sketch}: The Mellin transform of $\mathrm{Tr}\,
\mathrm{F}$ exists by assumption so
$e^{-\mathcal{M}_{\lambda} \left[\mathrm{Tr}\, \mathrm{F};\alpha\right]}$
represents an absolutely convergent series for $\alpha\in\mathbb{S}_\lambda$. Hence,
\begin{eqnarray}
e^{-\mathcal{M}_{\lambda} \left[\mathrm{Tr}\, \mathrm{F};\alpha\right]}
=\sum_{n=0}^\infty\frac{(-1)^n}{n!}\mathcal{M}_{\lambda}
\left[\mathrm{Tr}\, \mathrm{F};\alpha\right]^n
&=&\sum_{n=0}^\infty\frac{(-1)^n}{n!}\mathcal{M}_{\lambda}
\left[(\mathrm{Tr} \,\mathrm{F})^n;\alpha\right]\notag\\
&=&\mathcal{M}_{\lambda}
\left[\sum_{n=0}^\infty\frac{(-1)^n}{n!}(\mathrm{Tr}\,
\mathrm{F})^n;\alpha\right]\notag\\
&=&\mathcal{M}_{\lambda}\left[\mathrm{E}^{-\mathrm{Tr}\,
\mathrm{F}};\alpha\right]
\end{eqnarray}
where moving the power of $n$ into the functional Mellin transform in the first line follows from induction on Proposition $4.4$ \cite{LA1} because the multiplication represented by
$(\mathrm{Tr}\, \mathrm{F})^n$ is the $\ast$-convolution and
$\mathrm{Tr}\, \mathrm{F}(g)\in\C$, i.e. $\mathfrak{C}^\ast$ is commutative.  Equality between the first and second
lines follows from the absolute convergence of
$e^{-\mathcal{M}_{\lambda} \left[\mathrm{Tr}\, \mathrm{F};\alpha\right]}$,
the existence of the integral
$\mathcal{M}_{\lambda}\left[\mathrm{E}^{-\mathrm{Tr}\,
\mathrm{F}};\alpha\right]$ for the common domain $\alpha\in\mathbb{S}_\lambda$, and the fact that $e^{-\mathrm{tr}\,\mathrm{F}(g)}$ is analytic. It should be stressed that the equality in the lemma holds \emph{only} for $\alpha\in\mathbb{S}_\lambda$ properly
restricted. $\qED$

\begin{corollary}\label{key corollary}
Under the conditions of Lemma \emph{\ref{key lemma}}, replace
$\mathrm{Tr}\, \mathrm{F}$ with
$\mathrm{V}\in \mathrm{Mor}_C(G^\C,\mathfrak{C}^\ast)$ where now
 $\mathrm{V}(g)=V\rho(g)$ is self-adjoint, then the lemma together with \emph{Proposition $4.5$ \cite{LA1}} imply
\begin{equation}
\mathcal{M}_{\lambda}\left[\mathrm{E}^{-\mathrm{V}};\alpha\right]
=e^{-\mathcal{M}_{\lambda}\left[\mathrm{V};\alpha\right]}\,,\;\;\;\;\;\alpha\in\R\cap \mathbb{S}_\lambda\;.
\end{equation}
\end{corollary}

To finish the proof, put $\mathrm{F}\equiv\mathrm{E}^{-\mathrm{A}}$ in the lemma
and note that $\mathrm{F}(g)=e^{-\mathrm{A}(g)}$ so
\begin{equation}
\mathcal{M}_{\lambda}\left[\mathrm{E}^{-\mathrm{Tr}\,
\mathrm{E}^{-\mathrm{A}}};\alpha\right]=\int_{G^\C} e^{-\mathrm{tr}\left(e^{-\mathrm{A}(g)}\right)}\,\det g^\alpha\,\mathcal{D}_\lambda g=(\mathrm{Det}\,\mathrm{E}^{-\mathrm{A}})_\lambda^{-\alpha}
=\mathrm{det}\,(\mathrm{E}^{-\mathrm{A}})_\lambda^{-\alpha}
\end{equation}
and
\begin{equation}
\mathcal{M}_{\lambda} \left[\mathrm{Tr}\,\mathrm{E}^{-\mathrm{A}};\alpha\right]
=\int_{G^\C}\mathrm{tr}\left(e^{-\mathrm{A}(g)}\rho (g^\alpha)\right)\,\mathcal{D}_\lambda g
=(\mathrm{Tr}\,A)_\lambda^{-\alpha}
=\mathrm{tr}\,A_\lambda^{-\alpha}\;.
\end{equation}
Hence
\begin{equation}
(\mathrm{Det}\,\mathrm{E}^{-\mathrm{A}})_\lambda^{-\alpha}
= e^{-\mathrm{tr}\,A_\lambda^{-\alpha}}\;.
\end{equation} $\QED$

If $\mathrm{E}^{-\mathrm{A}}(g)$ is self-adjoint,  the corollary implies
\begin{equation}
 (\mathrm{E}^{-\mathrm{A}})^{-\alpha}_\lambda=\mathcal{M}_{\lambda}
\left[\mathrm{E}^{-\mathrm{E}^{-\mathrm{A}}};\alpha\right]=e^{-\mathcal{M}_{\lambda}
\left[\mathrm{E}^{-\mathrm{A}};\alpha\right]}=e^{-A_\lambda^{-\alpha}}\;,
\end{equation}
and with the help of Proposition $5.2$ \cite{LA1} the theorem can be rewritten as
 \begin{equation}
e^{-\mathrm{tr}\,A_\lambda^{-\alpha}}=\mathrm{det}
\,e^{-A_\lambda^{-\alpha}} \;,\;\;\;\;\;\alpha\in \R\cap\mathbb{S}_\lambda\;
\end{equation}
which is just the functional form of the standard relation with $M\equiv -A_\lambda^{-\alpha}$.

\subsection{Riemann zeros}
Let $\alpha_0$ denote a zero of the Riemann zeta function.
\begin{claim}
\begin{equation}
\zeta_{\mathrm{A}_{\lambda}}(\alpha_0)
=0\;\Rightarrow\;\Re{(\alpha_0)}=1/2\;, \;\;\;\;\;\alpha_0\in\langle0,1\rangle\subseteq\mathbb{S}_\lambda\;.
\end{equation}
\end{claim}
\vskip 1em
\emph{Argument}:
Use Theorem \ref{key theorem} with
$G^\C=\Omega\C^\times$ and $\mathrm{A}(g)=\mathrm{deg}(g)\,g$ as defined in subsection $2.1$. We can use  $\zeta_{A_{\mathcal{C}}}(\alpha)
=\mathcal{M}_{\mathcal{C}}\left[\mathrm{Tr}\,\mathrm{E}^{-\mathrm{A}};\alpha\right]
=\zeta(\alpha)$ to write
\begin{equation}\label{theorem}
e^{- \zeta_{\mathrm{A}_{\mathcal{C}}}(\alpha)} =e^{- \mathrm{tr}
A_{\mathcal{C}}^{-\alpha}} 
=(\mathrm{Det}\,\mathrm{E}^{-\mathrm{A}})_{\mathcal{C}}^{-\alpha}
=\mathrm{det}\,(\mathrm{E}^{-\mathrm{A}})_{\mathcal{C}}^{-\alpha}
\end{equation}
where the \emph{joint} fundamental strip has yet to be determined.\footnote{It is easy to verify the second equality of (\ref{theorem}) in \emph{finite} dimensions for $\alpha\in(0,\infty)/\{1\}$ taking
\begin{equation}
{A^{-\alpha}}(N):=\mathrm{diag}
\big\{{1^{-2\alpha}},{2^{-2\alpha}},\ldots,{j^{-2\alpha}},\ldots,
{N^{-2\alpha}}\big\}
=\mathrm{diag}
\big\{{1^{2}},{2^{2}},\ldots,{j^{2}},\ldots,
{N^{2}}\big\}^{-\alpha}\;.
\end{equation}
The pivotal point is weather our methods show the equality holds in the infinite-dimensional case as well.}

It is clear from the theorem that a joint fundamental region that includes the critical strip $\langle0,1\rangle$ along with the involution symmetry of the zeta zeros will lead to a relation between two functional determinants. From the previous subsection, the left-hand side of (\ref{theorem}) is already well-defined for $\alpha\in\langle0,\infty\rangle\backslash\{1\}$. So the first order of business is to determine the fundamental strip for the right-hand side.

According to Proposition $5.2$ \cite{LA1}, for some $\mathrm{V}\in \mathrm{Mor}_C(G^\C,\mathfrak{C}^\ast)$  such that $\mathrm{V}(g)=Vg$ is trace class and $0\notin\sigma(V_\lambda)$, the functional determinant is
\begin{equation}
(\mathrm{Det}\mathrm{V})_{\lambda}^{-\alpha}
=\det\,V_\lambda^{-\alpha}
=\mathcal{N}_\lambda(\alpha)\,\det
V^{-\alpha}
\,,\;\;\;\;\;\alpha\in\mathbb{S}_\lambda
\end{equation}
where $\mathcal{N}_\lambda(\alpha)$ is a $\lambda$-dependent
normalization,  and $\mathrm{det}\,V^{-\alpha}$ is well-defined for $\alpha\in\mathbb{S}_\lambda$ because $Vg$ is trace class. Since here $\mathrm{V}(g)\equiv\mathrm{E}^{-\mathrm{A}(g)}$ is self-adjoint, by Corollary $2.1$ this can be written more explicitly;
\begin{equation}
(\mathrm{Det}\,\mathrm{E}^{-\mathrm{A}})_\lambda^{-\alpha}
=\mathcal{N}_\lambda(\alpha)\,\det
e^{-A^{-\alpha}}
\,,\;\;\;\;\;\alpha\in\R\cap\mathbb{S}_\lambda\;.
\end{equation}

But $\mathcal{N}_\lambda(\alpha)$ is not well-defined since $Id$ is not trace class. So put $Id=RR^{-1}$ with $R$ some invertible, self-adjoint trace-class operator and invoke zeta function regularization to get (using results from \cite[\S\S4,5,7]{LA1})
\begin{eqnarray}
\mathcal{M}_{\mathcal{C}}[\mathrm{E}^{-\mathrm{Tr}\,\mathrm{Id}};\alpha]\,
=\mathrm{det}\,(\mathrm{Id})^{-\alpha}_{\mathcal{C}}
=\mathrm{det}\,(\mathrm{R}\ast\mathrm{R}^{-1})^{-\alpha}_{\mathcal{C}}
&=&\mathrm{det}\,(\mathrm{R})^{-\alpha}_{\mathcal{C}}\cdot
\mathrm{det}\,(\mathrm{R}^{-1})^{-\alpha}_{\mathcal{C}}\notag\\
&=&(\mathrm{det}\,{R}\cdot\mathrm{det}\,{R}^{-1})^{-\alpha}\notag\\
&=&(e^{-\zeta'_{\mathrm{R}_{\mathcal{C}}}(0)-\zeta'_{\mathrm{R}_{\mathcal{C}}^{-1}}(0)})^{-\alpha}\notag\\
&=&1^{-\alpha}
\end{eqnarray}
where we used $\log A:=2\sum_{m=3}^\infty(-1)^{m+1}C(A)^m=-\log A^{-1}$ with $C(A)=(A+Id)(A-Id)^{-1}$ (which is valid for $\sigma(A^{-1})\subseteq\C_+$\cite{HIG}) in the exponent. Hence,
\begin{equation}
\mathcal{N}_{\mathcal{C}}(\alpha)
=\int_{\mathcal{C}}e^{-\mathrm{tr}\,Id\,g}\,\mathrm{det}\,g^\alpha \;d\nu(g_{\mathcal{C}})
\equiv1\;.
\end{equation}

It follows that the theorem is valid
in the common fundamental strip
$\alpha\in\R\cap\langle0,\infty\rangle\backslash\{1\}$, and in particular it is \emph{single-valued} on the interval $\alpha\in\R\cap\langle0,1\rangle=(0,1)$;
\begin{equation}
e^{-\zeta_{\mathrm{A}_{\mathcal{C}}}(\alpha)} =\mathrm{det}\,e^{-{A^{-\alpha}}}\,\,\,\,\,\,\alpha\in(0,1)\;.
\end{equation}

Finally, being single-valued in the critical strip together with the involution symmetry of the zeros implies
\begin{equation}
\begin{array}{c}
  \zeta_{\mathrm{A}_{\mathcal{C}}}
(\alpha_0)=0=\zeta_{\mathrm{A}_{\mathcal{C}}}(1-\alpha^\ast_0)
\;,\;\;\;\;\;\alpha_0\in\langle0,1\rangle \\
  \big\Downarrow \\
  \mathrm{det}
\,e^{-A^{-\Re(\alpha_0)}} =
\mathrm{det}
\,e^{-A^{\Re(\alpha^\ast_0)-1}}\\
\big\Downarrow\\
 1=
\mathrm{det}\,e^{-A^{\Re(\alpha^\ast_0)-1}}/
\mathrm{det}\,e^{-A^{-\Re(\alpha_0)}}\;.
\end{array}
\end{equation}
Conclude that $\Re{(\alpha_0)}=1/2$ since $\mathrm{A}=\mathrm{A}^\ast$, there are no non-trivial zero modes of $\mathrm{A}$, and $\mathcal{N}_{\mathcal{C}}(\alpha)=\mathcal{N}_{\mathcal{C}}(1-\alpha)$ for $\alpha\in\langle0,1\rangle$.

Since diagonal $A(g)$ allowed for $\rho(g)\in Z(\mathfrak{C}^\ast)$, the theorem actually holds for complex $\alpha\in\langle0,\infty\rangle\backslash\{1\}$ and is single valued for complex $\alpha\in\langle0,1\rangle$. This doesn't change the conclusion. Also, the fundamental strip can't be analytically extended to the left of the imaginary
axis for both sides of (\ref{theorem}) simultaneously, i.e. there isn't a
\emph{consistent} choice of $\nu(g_\lambda)$ applicable to both
$\mathrm{Tr}\,\mathrm{A}$ and $\mathrm{E}^{-\mathrm{Tr}\,\mathrm{A}}$ to the
left of the imaginary axis. Furthermore, choosing a different normalization doesn't spoil the conclusion. For example the Haar measure $\nu(g_\lambda)=\log(g_\lambda)$ for $\R_+$ yields $\mathcal{N}_{\mathcal{C}}(\alpha)=\frac{2\pi i \Gamma(\alpha)}{\pi\csc(\alpha)}=2i\Gamma(1-\alpha)$. In the final step, this leads to an extra factor $\Gamma(1-\alpha_0)/\Gamma(\alpha_0)
=\Gamma(1-\alpha_0)\Gamma(\alpha_0^\ast)/|\Gamma(\alpha_0)|^2$ which also implies $\Re(\alpha_0)=1/2$.
$\QED$

Besides relying on the validity of the key theorem, expressing the functional determinant as a well-defined, regularized determinant is crucial. At the risk of repeating, if there is a disqualifying error in our arguments, it seems likely to be at this step.

The strategy used here can be applied to other zeta functions
endowed with appropriate reflection symmetry with obvious
modifications. In particular, for the Dirichlet eta function it easily follows that
$\eta_\mathrm{A}(\alpha_0)=0\Rightarrow\Re{(\alpha_0)}=1/2$.

\end{document}